\documentclass{article}
\usepackage{amssymb}

\input{tcilatex}
\begin{document}

\begin{center}
{\LARGE Range results for some social choice correspondences\medskip }

Jerry S. Kelly\footnote{%
Department of Economics, Syracuse University. \ Email: jskelly@syr.edu}%
{\LARGE \medskip }

December 16, 2017\bigskip 
\begin{eqnarray*}
&&\text{Abstract} \\
&&\text{Determination of the range of} \\
&&\text{several standard social choice} \\
&&\text{correspondences; Borda rule,} \\
&&\text{plurality rule, top-cycle, the} \\
&&\text{Copeland rule, approval} \\
&&\text{voting.}
\end{eqnarray*}%
\medskip 
\end{center}

Social choice correspondences might be partially evaluated in terms of what
they \textit{do}, in terms of what sets get chosen. Presumably we would not
think well of a rule that only selected either one specific singleton $\{a\}$%
, or the set $X$ of all alternatives. Range conditions could also be part of
characterization results.\ In this paper we determine the range of several
standard social choice correspondences.\medskip

The is a wide variety of such determinations. \ Assuming strong preference
orderings, any resolute, neutral correspondence has a range consisting of m
singleton sets. \ Many rules, like Pareto or Borda have a range consisting
of all non-empty subsets of alternatives. \ The Copeland Rule falls just
barely short of that full range property. \ In between, the range of
plurality voting depends on detailed relations between the number of
individuals and the number of alternatives.\medskip

Let $X$ with cardinality $|X|=m\geq 3$ be the set of \textbf{alternatives}
and let $N=\{1,2,...,n\}$ with $n\geq 2$ be the set of \textbf{individuals}.
A (strong) \textbf{ordering} on $X$ is a complete, asymmetric, transitive
relation on $X$ (non-trivial individual indifference is disallowed). The
highest ranked element of an ordering $r$ is denoted $r[1]$, the second
highest is denoted $r[2]$, etc. Also\ $r[1:k]$ is the set of alternatives in
the top $k$ ranks of $r$. \ The set of all orderings on $X$ is $L(X)$. A 
\textbf{profile} $u$ is an element $(u(1),u(2),...,u(n))$ of the Cartesian
product $L(X)^{N}$. \ If $x$ ranks above $y$ in $u(i)$, we write $x\succ
_{u(i)}y$.\medskip

A \textbf{social choice correspondence} $G$ is a map from the domain\ $%
L(X)^{N}$ to non-empty subsets of $X$. For a discussion of most of the
social choice correspondences in this paper, see, for example Schwartz
(1986) and Sen (2017). The \textbf{range} of social choice correspondence $G$
is the collection of all sets $S$ such that there exists a profile $u$ with $%
G(u)=S$.\medskip

Some range determinations are trivial.\medskip

\qquad (1) A constant social choice correspondence, say $G(u)=\{x,y\}$ for
all $u$ obviously has a singleton range.\medskip

\qquad (2) With strong orderings, a dictatorial social choice correspondence
has a range consisting of the $m$ different singleton sets.\medskip

\qquad (3) If $G(u)=T(u)=\cup _{i}u(i)[1]$, the set $T(u)$ of the tops of
all individual orderings at $u$, the range is all non-empty subsets of $X$
of size less than or equal to $n$.\medskip

For an example that is just above trivial, consider the Pareto optimal
correspondence. \ Alternative $y$ is \textbf{Pareto dominated} by another
alternative $x$ at profile $u$ if, for all $i$, individual $i$ ranks $x$
above $y$ in $u(i)$. \ Alternative $y$ is \textbf{Pareto optimal} at $u$ if
it is Pareto dominated by no other alternative at $u$. \ Consider $%
G_{Pa}(u)=\{x:x$ is Pareto optimal at $u\}$. \ Then $G_{Pa}(u)\neq
\varnothing $ since $T(u)\subseteq G(u)$. \ The range of the Pareto
correspondence consists of all non-empty subsets $S$ of $X$. \ Suppose $%
|S|=k $. \ To construct a profile $u$ with $G_{Pa}(u)=S$, start by putting
all of $X\backslash S$ in the bottom $m-k$ ranks in any order for everyone.
Then, for $\#1$, put $S$ in the top $k$ ranks in any ordering $O$. \ Follow
that by putting $S$ in the top $k$ ranks for $\#2$ but ordered as in $O^{-1}$%
. (Here $O^{-1}$ is the inverse of $O$: \ $xO^{-1}y$ if and only if $yOx$.)
\ Finally fill in the rest of the profile in any manner.\medskip

For another example at this level, we treat the maximin correspondence which
selects the alternatives whose worst ranking (over individuals) is highest.
\ Given profile $u$, let $n(u,x)$ be the largest integer $p$ for which there
is an individual $i$ with $u(i)[p]=x$. Then the \textbf{maximin
correspondence} sets%
\[
G_{M}(u)=\{x:n(u,x)\leq n(u,y)\text{ for all }y\in X\}\text{.} 
\]%
\newline
By neutrality, if we can find a set of $k$ elements in the range of $G_{M}$,
all sets of $k$ elements are in the range.\medskip

For $k=1$, we start with the profile $u_{1}$\medskip

\begin{tabular}{|l|l|l|l|l|}
\hline
$1$ & $2$ & $3$ & $\cdots $ & $n$ \\ \hline
$x$ & $x$ & $x$ &  & $x$ \\ 
$y$ & $y$ & $y$ & $\cdots $ & $y$ \\ 
$z$ & $z$ & $z$ &  & $z$ \\ 
$w$ & $w$ & $w$ & $\cdots $ & $w$ \\ 
$\vdots $ & $\vdots $ & $\vdots $ &  & $\vdots $ \\ \hline
\end{tabular}%
\medskip \newline
with $G_{M}(u_{1})=\{x\}$.\medskip

For $k=2$, we construct $u_{2}$ from $u_{1}$ by switching the alternatives
in first and second ranks for individual \#1:\medskip

\begin{tabular}{|l|l|l|l|l|}
\hline
$1$ & $2$ & $3$ & $\cdots $ & $n$ \\ \hline
$y$ & $x$ & $x$ &  & $x$ \\ 
$x$ & $y$ & $y$ & $\cdots $ & $y$ \\ 
$z$ & $z$ & $z$ &  & $z$ \\ 
$w$ & $w$ & $w$ & $\cdots $ & $w$ \\ 
$\vdots $ & $\vdots $ & $\vdots $ &  & $\vdots $ \\ \hline
\end{tabular}%
\medskip \newline
with $G_{M}(u_{2})=\{x,y\}$.

For $k=3$, we construct $u_{3}$ from $u_{2}$ by switching the alternatives
in second and third ranks for individuals \#1 and \#2:\medskip

\begin{tabular}{|l|l|l|l|l|}
\hline
$1$ & $2$ & $3$ & $\cdots $ & $n$ \\ \hline
$y$ & $x$ & $x$ &  & $x$ \\ 
$z$ & $z$ & $y$ & $\cdots $ & $y$ \\ 
$x$ & $y$ & $z$ &  & $z$ \\ 
$w$ & $w$ & $w$ & $\cdots $ & $w$ \\ 
$\vdots $ & $\vdots $ & $\vdots $ &  & $\vdots $ \\ \hline
\end{tabular}%
\medskip \newline
with $G_{M}(u_{3})=\{x,y,z\}$. Continue in this fashion until $k$ is the
smaller of $m$ and $n$.\medskip

Now we take up more complicated cases.\medskip

\section{The Borda correspondence\protect\medskip}

At profile $u$, let $s(x,i,u)$ be the rank of alternative $x$ in $i$'s
ordering (so if, at $u$, individual $i$ has $x$ top-ranked, $s(x,i,u)=1$). \
The Borda score for $x$ at $u$ is the sum\medskip

\[
B(x,u)=\Sigma _{i}s(x,i,u) 
\]%
\newline
of those ranks over all individuals. \ Then the \textbf{Borda correspondence}
selects the alternatives with minimal Borda score:\medskip

\[
G_{B}(u)=\{x\in X:B(x,u)\leq B(y,u)\text{ for all }y\in X\} 
\]%
\bigskip

If $n$ is even, the range of Borda consists of all possible subsets $S$ of $%
X $. \ To get $S$, construct the profile $u$ as follows: Let $P$ be an
ordering of $S$ and $Q$ be an ordering of $X$. $u(1)$ will be $P$ above and $%
X\backslash S$ below, ordered in any way. \ $u(2)$ will be $P^{-1}$ above
and $X\backslash S$ below, ordered in any way. \ Then half of the remaining
individuals have ordering $Q$ while the other half have ordering $Q^{-1}$. \
If $n$ is odd, there is a profile $u$ with $G_{B}(u)=S$ unless $S=X$ where $%
|X|$ is even. \ For proof of this, see Kelly and Qi (2016), where also
partial results are given for the range of the \textbf{Borda social welfare
function}, which yields social \textit{rankings}.\medskip

\section{Plurality rule\protect\medskip}

Given a profile $u$ and an element $x$ of $X$, let $n_{x}(u)$ be the number
of individuals with $x$ top-ranked at $u$. Alternative $x$ is a \textbf{%
plurality winner} at $u$ if $n_{x}(u)\geq n_{z}(u)$ for all $z$ in $X$. 
\textbf{Plurality rule} is the social choice correspondence $G_{Pl}$ that,
at each profile $u$, has $G_{Pl}(u)$ equal to the set of all plurality
winners at $u$.\ \ The \textbf{plurality number} at $u$, $N(u)$, is the
value of $n_{x}$ for a plurality winner $x$ at $u$. \ When profiles $u$ and $%
u^{\ast }$ are fixed for a discussion, we simplify $n_{x}(u)$ and $%
n_{x}(u^{\ast })$ respectively to $n_{x}$ and $n_{x}^{\ast }$.\medskip\ 

For the range of plurality rule, let $S$ be a subset of $X$ so $|S|$ $=k\leq
m$. \ When is $S=\{x_{1},x_{2},...,x_{k}\}$ in the range of plurality rule?
We certainly need $k\leq n$ since there are at most $n$ tops. \ If $k|n$,
let $u$ be a profile where the first $n/k$ individuals have $x_{1}$ on top,
the next $n/k$ have $x_{2}$ on top, and so on. \ $S=G_{Pl}(u)$. \ If $k\nmid
n$, let $n=kq+r$ with $0<r<k$. \ Let \ $u$ be a profile where the first $q$
individuals have $x_{1}$ on top, the next $q$ have $x_{2}$ on top, and so
on. If $k=m$ and $r>0$, sets $S$ with $|S|=k$ are not in the range.\medskip

If $m>k$, there are $r$ individuals whose tops at $u$ are still undetermined
and $m-k$ alternatives to use for those tops. \ For $x_{1}$, $x_{2}$, ..., $%
x_{k}$ to be the plurality winners, we would need $q$ strictly larger than
the highest number needed for a non-winner:

\[
q>\left\lceil \frac{r}{m-k}\right\rceil =\left\lceil \frac{n-qk}{m-k}%
\right\rceil 
\]%
\newline
where $\left\lceil z\right\rceil $ is the smallest integer not less than $z$%
. \ Here $q$ is $\left\lfloor \frac{n}{k}\right\rfloor $ where $\left\lfloor
z\right\rfloor $ is the largest integer not greater than $z$. \ Combining,
when

\[
\left\lfloor \frac{n}{k}\right\rfloor >\left\lceil \frac{n-\left\lfloor 
\frac{n}{k}\right\rfloor k}{m-k}\right\rceil 
\]%
\newline
any set of $k$ elements is in the range of plurality rule.\medskip

\section{The top cycle correspondence\protect\medskip}

Here $G_{T}(u)$ is the set of maximal alternatives for the \textit{%
transitive closure} of the simple majority voting relation.\ Given $x,y\in X$%
, we define the \textbf{simple majority voting} (SMV) relation as $x\succsim
_{u}y$ just when $|\{i:x\succ _{u(i)}y\}|$ $\geq |\{i:y\succ _{u(i)}x\}|$. \
Let $R(u)$ be the \textbf{transitive closure} of $\succsim _{u}$. \ That is, 
$xR(u)y$ if there is a sequence $x_{1}$, $x_{2}$, ..., $x_{t}$ such that%
\[
x\succsim _{u}x_{1}\succsim _{u}x_{2}\succsim _{u}...\succsim
_{u}x_{t}\succsim _{u}y 
\]%
\newline
If $xR(u)y$ and $yR(u)x$, we write $x\approx y$. Then $G_{T}(u)=\{x:xR(u)y$
for all $y\}$. \ Note that if $y\in G_{T}(u)$ and $xR(u)y$, then $x\in
G_{T}(u)$.\medskip

Note first, this rule is neutral, so if we can find a set of $k$ elements in
the range, all sets of $k$ elements are in the range.\medskip

\textit{Case 1}. $n$ odd. For odd $n$, any one element set is in the range,
but no two element set is. Now we construct, for any subset $S$ of $X$
consisting of three or more alternatives, a profile $u$ that has $S$ in the
range. \ First we observe that with three individuals, there is a profile on
just $S$ such that the transitive closure of SMV is social indifference on $%
S $.\ The standard voting paradox shows a profile $v$ on three individuals
and exactly three alternatives that exhibits a voting cycle: $x\succ y\succ
z\succ x$, and the transitive closure is $x\approx y\approx z\approx x$. \
For $k=4$, pick one alternative from $\{x,y,z\}$, say $y$, and a new
alternative, $w$. Now insert $w$ in each ordering in profile $v$ just below $%
y$. That creates a 4-element cycle and the transitive closure of SMV is
social indifference among these 4 alternatives. Continuing in this manner,
we can construct a profile $s$ such that the transitive closure of SMV is
social indifference on $S$. \ Now let $P$ be an arbitrary ordering of $%
X\backslash S$ and $Q$ is an arbitrary ordering of $X$, and construct $u$ as
follows:\medskip

\qquad A. \ For individuals \#1-\#3, put the elements of $S$ in the top $k$
ranks, ordered as in the profile $s$ on $S$ and put the elements of $%
X\backslash S$ in the bottom $m-k$ ranks ordered as in $P$;\medskip

\qquad B. \ For odd $j$, with $3<j\leq n$, set $u(j)=Q$;\medskip

\qquad C. For even $j$, with $2<j\leq n-1$, set $u(j)=Q^{-1}$.\medskip

Then $S$ is the maximal set of the transitive closure of SMV at $u$.\medskip

\qquad \qquad 
\begin{tabular}{|l|l|l|}
\hline
1 & 2 & 3 \\ \hline
$x$ & $z$ & $y$ \\ 
$y$ & $x$ & $w$ \\ 
$w$ & $y$ & $v$ \\ 
$v$ & $w$ & $z$ \\ 
$z$ & $v$ & $x$ \\ 
$\vdots $ & $\vdots $ & $\vdots $ \\ \hline
\end{tabular}%
\medskip \newline
where $G_{T}(u)=\{x,y,z,w,v\}$. This example also illustrates that $G$ does 
\textbf{not} satisfy the Pareto condition.\medskip

\textit{Case 2}. \ $n$ even. \ For even $n$, every subset is in the range. \
Let $S$ be any $k$-element subset of $X$ for $0<k\leq m$. \ Construct
profile $u$ as follows, where $O$ is an arbitrary ordering of $S$, $P$ is an
arbitrary ordering of $X\backslash S$ and $Q$ is an arbitrary ordering of $X$%
:\medskip

\qquad A. For \#1, put the elements of $S$ in the top $k$ ranks, ordered as
in $O$ and put the elements of $X\backslash S$ in the bottom $m-k$ ranks
ordered as in $P$;\medskip\ \ 

\qquad B. \ For \#2, put the elements of $S$ in the top $k$ ranks, ordered
as in $O^{-1}$ and put the elements of $X\backslash S$ in the bottom $m-k$
ranks ordered as in $P$;\medskip

\qquad C. \ For odd $j$, with $1<j\leq n-1$, set $u(j)=Q$;\medskip

\qquad D. For even $j$, with $2<j\leq n$, set $u(j)=Q^{-1}$.\medskip

Then $S$ is the maximal set of the transitive closure of SMV at $u$.\medskip

\section{The Copeland rule\protect\medskip}

Given a profile $u$ and an element $x$ of $X$, let $n_{x}(u)$ be the \textbf{%
Copeland score}$\footnote{%
There is a variant of the Copeland score that is the number of alternatives
defeated by $x$ \textit{minus} the number of alternatives that defeat $x$. \
The distinction between these variants is of no consequence for our range
results.}$, the number of alternatives defeated by $x$ under simple majority
vote at $u$. Alternative $x$ is a \textbf{Copeland winner} at $u$ if $%
n_{x}(u)\geq n_{z}(u)$ for all $z$ in $X$. The\textbf{\ Copeland rule} is
the social choice correspondence $G_{C}$ that, at each profile $u$, has $%
G_{C}(u)$ equal to the set of all Copeland winners at $u$.\medskip

Our focus is on the range of $G_{C}$. \ The Copeland correspondence is
neutral and by neutrality, a set $S$ from $X$ is in the range if and only if
there is at least one set $T\subseteq X$ in the range of $G_{C}$ with $|S|$ $%
=|T|$. \ Thus, by an abuse of language, we will say that $k$ is in the range
of $G_{C}$ if there is a $u$ with $|G_{C}(u)|$ $=k$. \ If $k$ is in the
range for $m$ and $n$, it is also in the range for $m+1$ and $n$ (just add
an alternative at the bottom of everyone's ordering). \ Also, if $k$ is in
the range for $m$ and $n$, it is also in the range for $m$ and $n+2$ (just
add two individuals with inverse orderings).\medskip

\textbf{Even n }\medskip

With even $n$, SMV ties are possible. \ That makes the determination of the
range much easier. \ We want to show for any $m$, $n$, and \textbf{any} set $%
S$ of $k$ alternatives where $1\leq k\leq m$, there is a profile at which
the Copeland choice set is $S$. \ By the first remarks in this section, it
suffices to show this for $n=2$. \ Let $P$ be any strong ordering on $S$ and 
$Q$ be any strong ordering on $X\backslash S$. \ Construct $u$ by setting
\#1's ordering to be $P$ on top followed by $Q$ and \#2's ordering to be $%
P^{-1}$ on top followed by $Q$.\medskip

\textbf{Odd n}\medskip

Now assume all preferences are strong orderings on $X$ and there are an odd
number of individuals. \ Given profile $u$, consider the $m$-term
non-increasing sequence $(s_{1},s_{2},...s_{m})$ obtained from $n_{x}(u)$, $%
n_{y}(u)$, ... by rearranging according to magnitude, largest to
smallest.\medskip

The set $G_{C}(u)$ of maximal elements of $R(u)$ will have $k$ elements just
when the associated sequence $s_{1},s_{2},...s_{m}$ has an initial
subsequence of $k$-elements, all the same, and then, if $k<m$, a $k+1^{st}$
element different from the first $k$. \ Note that since there are $m(m-1)/2$
pairs of distinct alternatives, $s_{1}+s_{2}+...+s_{m}=m(m-1)/2$.\medskip

For $|X|=m$, there exists at least one $k$, with $1\leq k\leq m$ such that $%
k $ is \textbf{not} in the range of $G_{C}$. \ To be more precise\medskip

\qquad (1) \ If $m$ is even, $m$ is not in the range of $G_{C}$;

\qquad (2) \ If $m$ is odd, $m-1$ is not in the range of $G_{C}$.\medskip

For (1), note that for $m$ to be in the range of $G_{C}$, we would have to
have a profile $u$ with associated sequence $%
(s_{1},s_{2},...s_{m})=(a,a,...,a)$ for some positive integer $a$. \ But
then we would have $s_{1}+s_{2}+...+s_{m}=ma=m(m-1)/2$ so $a=(m-1)/2$ which
is only possible if $m$ is odd.\medskip

For (2), note that for $m-1$ to be in the range of $G_{C}$, we would have to
have a profile $u$ with associated sequence $%
(s_{1},s_{2},...s_{m})=(a,a,...,a,b)$ for some positive integer $a,b$ with $%
b<a$. \ Then $s_{1}+s_{2}+...+s_{m}=(m-1)a+b=m(m-1)/2$ with $0\leq b<a$. \
Thus $(m-1)a\leq m(m-1)/2<ma$ and so $(m-1)/2<a\leq m/2$, which is not
possible for an integer $a$ if $m$ is odd.\medskip

We now show that (1) and (2) are the only exceptions: For $n=3$ (and so also
for all larger odd $n$), and even $m>2$, every $k$, $1\leq k<m$ is in the
range of $G_{C}$. \ For $n=3$ (and so also for all larger odd $n$), and odd $%
m$, every $k$, $1\leq k\leq m$ except $m-1$, is in the range of $G_{C}$%
.\medskip

By what we have already argued, it suffices to show that $k$ appears in the
range "as soon as possible," i.e., for odd $m$, we have $k=m$ in the range
for $m$ alternatives and for even $m$, we have $k=m$ in the range for $m+2$
alternatives.\medskip

\textbf{Part 1: m odd.} \ To get $k=m$ in the range for odd $m$, we
construct profiles where the Copeland score for all $m$ alternatives is $1$.
\ First, note that $k=3$ is in the range for $m=3$ at the profile $u_{1}$%
:\medskip

\qquad \qquad 
\begin{tabular}{|l|l|l|}
\hline
1 & 2 & 3 \\ \hline
$x$ & $y$ & $z$ \\ 
$y$ & $\mathbf{z}$ & $x$ \\ 
$z$ & $\mathbf{x}$ & $\mathbf{y}$ \\ \hline
\end{tabular}%
\medskip \newline
(the classic voting paradox). \ For later comparison, we have emboldened the
bottom two alternatives in $u_{1}(2)$ and the bottom alternative in $%
u_{1}(3) $. These emboldened alternatives make up the whole set of
alternatives with no duplicates.\medskip

For $m=5$, we construct $u_{2}$ from $u_{1}$ by inserting two new
alternatives, $a$ and $b$. \ Alternative\ $a$ is inserted at the top of 1's
ordering and at the bottom - emboldened - for \#3. \ $b$ is inserted at the
bottom of \#1's ordering and at the top for \#2:\medskip

\qquad \qquad 
\begin{tabular}{|l|l|l|}
\hline
1 & 2 & 3 \\ \hline
$a$ & $b$ & $z$ \\ 
$x$ & $y$ & $x$ \\ 
$y$ &  &  \\ 
$z$ & $\mathbf{z}$ & $\mathbf{y}$ \\ 
$b$ & $\mathbf{x}$ & $\mathbf{a}$ \\ \hline
\end{tabular}%
\medskip \newline
Finally, $a$ is inserted at the middle position for \#2 and b is inserted in
the middle position for \#3 (that $b$ is emboldened):\medskip

\qquad \qquad 
\begin{tabular}{|l|l|l|}
\hline
1 & 2 & 3 \\ \hline
$a$ & $b$ & $z$ \\ 
$x$ & $y$ & $x$ \\ 
$y$ & $a$ & $\mathbf{b}$ \\ 
$z$ & $\mathbf{z}$ & $\mathbf{y}$ \\ 
$b$ & $\mathbf{x}$ & $\mathbf{a}$ \\ \hline
\end{tabular}%
\medskip \newline
The emboldened alternatives make up the whole set of alternatives with no
duplicates. \ Each alternative has Copeland score of 2.\bigskip

For $m=7$, we construct $u_{3}$ from $u_{2}$ by inserting two new
alternatives, $c$ and $d$. \ $c$ is inserted at the top of 1's ordering and
at the bottom - emboldened - for \#2. \ Alternative\ $d$ is inserted at the
bottom of \#1's ordering and at the top for \#3:\medskip

\qquad \qquad 
\begin{tabular}{|l|l|l|}
\hline
1 & 2 & 3 \\ \hline
$c$ & $b$ & $d$ \\ 
$a$ & $y$ & $z$ \\ 
$x$ & $a$ & $x$ \\ 
$y$ &  &  \\ 
$z$ & $\mathbf{z}$ & $\mathbf{b}$ \\ 
$b$ & $\mathbf{x}$ & $\mathbf{y}$ \\ 
$d$ & $\mathbf{c}$ & $\mathbf{a}$ \\ \hline
\end{tabular}%
\medskip \newline

Finally, $d$ and $c$ are inserted at the middle position for \#2 (where $d$
is emboldened) and \#3:\medskip

\qquad \qquad 
\begin{tabular}{|l|l|l|}
\hline
1 & 2 & 3 \\ \hline
$c$ & $b$ & $d$ \\ 
$a$ & $y$ & $z$ \\ 
$x$ & $a$ & $x$ \\ 
$y$ & $\mathbf{d}$ & $c$ \\ 
$z$ & $\mathbf{z}$ & $\mathbf{b}$ \\ 
$b$ & $\mathbf{x}$ & $\mathbf{y}$ \\ 
$d$ & $\mathbf{c}$ & $\mathbf{a}$ \\ \hline
\end{tabular}%
\medskip \newline
The emboldened alternatives make up the whole set of alternatives with no
duplicates. \ Each alternative has Copeland score of 3.\medskip

We show one more iteration. \ For $m=9$, we construct $u_{4}$ from $u_{3}$
by inserting two new alternatives, $r$ and $s$. \ $r$ is inserted at the top
of 1's ordering and at the bottom - emboldened - for \#3. \ Alternative\ $s$
is inserted at the bottom of \#1's ordering and at the top for \#2:\medskip

\qquad \qquad 
\begin{tabular}{|l|l|l|}
\hline
1 & 2 & 3 \\ \hline
$r$ & $s$ & $d$ \\ 
$c$ & $b$ & $z$ \\ 
$a$ & $y$ & $x$ \\ 
$x$ & $a$ & $c$ \\ 
$y$ &  &  \\ 
$z$ & $\mathbf{d}$ & $\mathbf{b}$ \\ 
$b$ & $\mathbf{z}$ & $\mathbf{y}$ \\ 
$d$ & $\mathbf{x}$ & $\mathbf{a}$ \\ 
$s$ & $\mathbf{c}$ & $\mathbf{r}$ \\ \hline
\end{tabular}%
\medskip \newline
Finally, $r$ and $s$ are inserted at the middle position for \#2 and \#3
(where $s$ is emboldened):\medskip

\qquad \qquad 
\begin{tabular}{|l|l|l|}
\hline
1 & 2 & 3 \\ \hline
$r$ & $s$ & $d$ \\ 
$c$ & $b$ & $z$ \\ 
$a$ & $y$ & $x$ \\ 
$x$ & $a$ & $c$ \\ 
$y$ & $r$ & $\mathbf{s}$ \\ 
$z$ & $\mathbf{d}$ & $\mathbf{b}$ \\ 
$b$ & $\mathbf{z}$ & $\mathbf{y}$ \\ 
$d$ & $\mathbf{x}$ & $\mathbf{a}$ \\ 
$s$ & $\mathbf{c}$ & $\mathbf{r}$ \\ \hline
\end{tabular}%
\medskip \newline
The emboldened alternatives make up the whole set of alternatives with no
duplicates. \ Each alternative has Copeland score of 4. \ This process,
alternating which of \#2 and \#3 has a new alternative inserted at the
bottom, can be continued indefinitely. \ For general odd $m$, at the
resulting profile every alternative has Copeland score $(m-1)/2$.\medskip

\textbf{Part 2: m even.} \ For even $m$, $k=m-2$ is in the range for $m$
alternatives. \ This requires only a very simple construction. \ Suppose we
use Part 1 to start with a profile $u$ for $m-1$ alternatives (making up set 
$X^{\ast }$) such that every alternative in $X^{\ast }$ has Copeland score $%
(m-2)/2$. \ \ Then construct profile $u^{\ast }$ by inserting a new
alternative, say $t$, as follows:

\qquad 1. \ $t$ is made the top-most alternative in $u^{\ast }(1)$;

\qquad 2. \ $t$ is made the bottom-most alternative in $u^{\ast }(2)$;$\geq $

\qquad 3. \ For individual \#3, $t$ is inserted into the next-to-bottom
space, just above the bottom element, call it $z$, in $u^{\ast }(3)$.\medskip

Then the Copeland scores of the $m-2$ alternatives in $X^{\ast }$ other than 
$z$ are increased by $1$ since they all now defeat $t$. \ The Copeland score
of $z$ is unchanged. \ Since $t$ only defeats $z$, the Copeland score of $t$
is $1$ less than the maximal Copeland score. \ The image of $G_{C}$ at $%
u^{\ast }$ is $X^{\ast }\backslash \{z\}$, of cardinality $m-2$. \medskip

\section{Approval voting\protect\medskip}

Approval voting (Brams and Fishburn, 1982) is not a collective choice
correspondence since different information is used. \ For a fixed $B$, where 
$1\leq B\leq m$, an \textbf{index vector} is an element of $B^{N}$. \ An 
\textbf{extended social choice correspondence} is a correspondence $%
G:L(X)^{N}\times B^{N}\rightarrow X$. \ Given preference profile $%
u=(u_{1},...,u_{n})$ and index vector $b=(b_{1},...,b_{n})$, the \textbf{%
approval voting score} for alternative x is $A_{S}(x,u,b)=|\{i:x\in
u(i)[1:b_{i}]\}|$, the number of individuals $i$ who have $x$ in their top $%
b_{i}$ ranks. \ \textbf{Approval voting} then is the extended social choice
correspondence $G_{A}$ that sets $G_{A}(u,b)=\{x:A_{S}(x,u,b)\geq
A_{S}(y,u,b)$ for all $y$ in $X\}$ where $B=M$, i.e., individuals can
approve any number of alternatives.\medskip

The range of $G_{A}$ is all non-empty subsets of $X$. \ Given subset $%
S\subseteq X$, with $|S|=k\geq 1$, let $u$ be a profile where every
individual has $S$ in their top $k$ ranks (in any order) and $b=(k,k,...,k)$%
. \ Then $G_{A}(u,b)=S$.\medskip

But one can generally get $G_{A}(u,b)=S$ for index vectors much smaller than 
$(k,k,...,k)$. \ Given index vector $b$, let the gauge of $b$ be the maximal
value of the $b_{i}$ components.\medskip

\textbf{Proposition.} \ Let $|S|=k\geq 1$.\medskip

\qquad (1) \ If $k\geq n$, there exists a profile $u$ and an index vector $b$
of gauge $g=\lceil k/n\rceil $ with $G_{A}(u,b)=S$;

\qquad (2) \ If $k<n$, there exists a profile $u$ and an index vector $b$ of
gauge $g\leq 2$ with $G_{A}(u,b)=S$.\medskip

\textbf{Proof}: \ For (1), let $u$ be a profile where the $k$ elements of $S$
are strung out in the top $\lceil k/n\rceil $ ranks, with one occurrence of
each (here illustrated with $k=2n+2$):\medskip

\qquad \qquad 
\begin{tabular}{|l|l|l|l|l|l|}
\hline
$1$ & $2$ & $3$ & $4$ & $\cdots $ & $n$ \\ \hline
$x_{1}$ & $x_{2}$ & $x_{3}$ & $x_{4}$ &  & $x_{n}$ \\ 
$x_{n+1}$ & $x_{n+2}$ & $x_{n+3}$ & $x_{n+4}$ & $\cdots $ & $x_{2n}$ \\ 
$x_{2n+1}$ & $x_{2n+2}$ & $\vdots $ & $\vdots $ &  & $\vdots $ \\ 
$\vdots $ & $\vdots $ &  &  &  &  \\ \hline
\end{tabular}%
\medskip \newline
and set $b=(\lceil k/n\rceil ,\lceil k/n\rceil ,\lceil k/n\rceil
-1,...,\lceil k/n\rceil -1)$.\medskip

\qquad For (2), where $n=kq+r$ for $0\leq r<n$ some care must be taken with
assumptions about $m$.\medskip

\textbf{Example. }\ Suppose $n=10$ and $k=4$. \ Here $q=2$ and you might
consider setting $u$ to be\medskip

\qquad \qquad 
\begin{tabular}{|l|l|l|l|l|l|l|l|l|}
\hline
$1$ & $2$ & $3$ & $4$ & $\cdots $ & $7$ & $8$ & $9$ & $10$ \\ \hline
$x_{1}$ & $x_{1}$ & $x_{2}$ & $x_{2}$ &  & $x_{4}$ & $x_{4}$ & $x_{5}$ & $%
x_{6}$ \\ 
$\vdots $ & $\vdots $ & $\vdots $ & $\vdots $ &  &  &  &  & $\vdots $ \\ 
\hline
\end{tabular}%
\medskip \newline
with $b=(1,1,...,1)$. \ The problem with this is that while we know $m\geq 4$
(since $k=4$), we don't know $m\geq 6$.\medskip

\qquad We deal with this by enlarging the number of occurrences of an
element of $S$ from $q$ to $q+1$ and spilling over to the second
rank:\medskip

\qquad \qquad 
\begin{tabular}{|l|l|l|l|l|l|l|l|l|}
\hline
$1$ & $2$ & $3$ & $4$ & $\cdots $ & $7$ & $8$ & $9$ & $10$ \\ \hline
$x_{1}$ & $x_{1}$ & $x_{1}$ & $x_{2}$ &  & $x_{3}$ & $x_{3}$ & $x_{3}$ & $%
x_{4}$ \\ 
$x_{4}$ & $x_{4}$ & $\vdots $ & $\vdots $ &  &  &  & $\vdots $ & $\vdots $
\\ 
$\vdots $ & $\vdots $ &  &  &  &  &  &  &  \\ \hline
\end{tabular}%
\medskip \newline
with $b=(2,2,1,...,1)$.\medskip

Generalizing, consider first the case of $k<n/2$. \ Thus $2k<n$ and in $%
n=kq+r$ we have $q\geq 2$ and $r<k$. \ We get bounds on k(q+1).\medskip

\qquad \qquad (i) \ $n=kq+r<kq+k=k(q+1)$;

\qquad \qquad (ii) \ $2n=2kq+2r>k(q+1)+2r>k(q+1)$.\medskip

Thus we can string $q+1$ occurrences of each of $k$ alternatives along the
top rank and then extend part way along the second rank and set $%
b=(2,2,...,2,1,...,1)$ as in the way we dealt with the example just
above.\medskip

Finally, we treat the case where $n/2<k<n$ and $q=1$. \ We again establish
bounds:\medskip

\qquad \qquad (i) \ $2k>n$ since $n/2<k$;

\qquad \qquad (ii) \ But $2k<2n$ since $k<n$.\medskip

Thus we can string $2=q+1$ occurrences of each of $k$ alternatives along the
top rank and then extend part way along the second rank and set $%
b=(2,2,...,2,1,...,1)$ \ \ \ \ \ $\square $

\bigskip \bigskip

\textbf{REFERENCES}\medskip

Brams, SJ and PK Fishburn (1982) \textit{Approval Voting} (Birkh\"{a}%
user).\medskip

Copeland, AH (1951) "A 'Reasonable' Social Welfare Function." (Univ. of
Michigan mimeo).\medskip

Kelly, JS and S Qi (2016) "A Conjecture on the Construction of Orderings by
Borda's Rule," \textit{Social Choice and Welfare} 47 113-126.\medskip

Schwartz, T (1986) \textit{The Logic of Collective Choice}
(Columbia).\medskip

Sen, A \ (2017) \textit{Collective Choice and Social Welfare} 2nd ed.
(Harvard).

\end{document}